\documentclass[12 pt]{article}

\pdfoutput=1

\usepackage[T1]{fontenc}
\usepackage[utf8]{inputenc}
\usepackage[english]{babel}
\usepackage{lmodern}
\usepackage{xspace} 

\title{\textbf{Dimension improvement in Dhar's refutation of the Eden conjecture}}
\author{Quentin \textsc{Bertrand} \& Jules \textsc{Pertinand}}
\date{May 2017}

\usepackage[numbers]{natbib}
\usepackage{graphicx}
\usepackage{color}
\usepackage{amssymb}
\usepackage{mathtools}
\usepackage[]{algorithm2e}
\usepackage{tikz}
\usepackage{stmaryrd}
\usepackage{hyperref}
\usepackage{amsthm}
\theoremstyle{definition}
\newtheorem{theorem}{Theorem}

\newcommand\card[1]{\left\vert{#1}\right\vert}

\DeclareFontFamily{U}{wncyr}{}
\DeclareFontShape{U}{wncyr}{m}{n}{<->wncyr10}{}
\DeclareFontShape{U}{wncyr}{m}{it}{<->wncyi10}{}
\DeclareFontShape{U}{wncyr}{m}{sc}{<->wncysc10}{}
\DeclareFontShape{U}{wncyr}{b}{n}{<->wncyb10}{}
\DeclareTextCommand{\guillemotleft}{T1}{%
  {\fontencoding{U}\fontfamily{wncyr}\selectfont\symbol{"3C}}%
}
\DeclareTextCommand{\guillemotright}{T1}{%
  {\fontencoding{U}\fontfamily{wncyr}\selectfont\symbol{"3E}}%
}

\usepackage[hmarginratio=1:1]{geometry}
\usepackage{enumitem}

\usepackage{float}

\usepackage{natbib}
\usepackage{graphicx}
\usepackage{color}
\usepackage{amssymb,amsmath,amsthm}
\usepackage{mathtools}
\usepackage{parskip}
\theoremstyle{definition}

\theoremstyle{remark}

\makeatletter
\def\thm@space@setup{%
  \thm@preskip= 0.3 cm \thm@postskip=0.2 cm
}
\makeatother

\begin{document}
 
\maketitle

\begin{abstract}
We consider the Eden model on the $d$-dimensional hypercubical unoriented lattice, for large $d$. 
Initially, every lattice point is healthy, except the origin which is infected. Then, each infected lattice point contaminates any of its neighbours with rate $1$. The Eden model is equivalent to first passage percolation, with exponential passage times on edges. The Eden conjecture states that the limit shape of the Eden model is a Euclidean ball.

By pushing the computations of Dhar \citep{Dhar} a little further with modern computers and efficient implementation we obtain improved bounds for the speed of infection. This shows that the Eden conjecture does not hold in dimension superior to $22$ (the lowest known dimension was $35$).

\end{abstract}


\section{The Eden model: definitions and previous results}

We consider the first passage percolation on a $d$-dimensional hypercubical unoriented lattice (\citep{refPerco}) as stated in \citep{Gerin}. Let $\{ \alpha(x,y) | (x, y) \in \text{edges of } \mathbb{Z}^d \}$ be a family of i.i.d random variables, with exponential law of parameter $1$. Let n $\in \mathbb{N}$.
For a path $\mathcal{W}$ : $x_0 \rightarrow  x_1 \rightarrow ... \rightarrow x_n $ of neighbouring vertices, we define the passage time along $\mathcal{W}$: $\alpha(\mathcal{W})=  \sum \limits_{i=1}^{n} \alpha(x_{i-1}, x_i) $. The family $\{ \alpha(x,y) | (x, y) \in \text{edges of } \mathbb{Z}^d \}$ defines a random distance, $\forall (x,y) \in \mathbb{Z}^d \times \mathbb{Z}^d $
$$ \mathcal{D}(x,y)=  \underset{ \mathcal{\mathcal{W}} \text{ path from } x \text{ to } y }{ \operatorname{inf}} \alpha( \mathcal{W} ).$$

For all $t \in \mathbb{R}$ we set $B_t = \{ x \in \mathbb{Z}^d | \mathcal{D}(0,x) \leq t\}$. Richardson (1973) and Cox-Durett (1981) have shown that there exists a compact convex $B^* \subset \mathbb{R}^d$  such that for all $\epsilon > 0 $ $$\mathbb{P}\left((1-\epsilon)B^* \subset \frac{B_t^d}{t} \subset (1+ \epsilon)B^*, \text{ for t big enough} \right)=1.$$ Eden conjectured that this limit form $B^*$ was a Euclidean ball in every dimension.

For all  $ n \in \mathbb{N} $ we note $\mathcal{P}_n^d$ the hyperplane of equation $x_1= n$ in dimension $d$. We note $\mathbf{0}$ the origin of the hypercube.
Observe that $\mathcal{D}(\mathbf{0},\mathcal{P}_n^d)$ is the distance between the origin and the hyperplan $\mathcal{P}_n^d$.
Cox-Durrett, Hammersley and Welsh (\citep{hammersley1965first}) have shown that with probability 1 $\lim\limits_{n \rightarrow +\infty} \dfrac{\mathcal{D}(\mathbf{0},\mathcal{P}_n^d)}{n}= \mu_{axis}^d $, for a certain $\mu_{axis}^d$, and moreover that 
$$\mu_{axis}^d = { \inf \limits_{n \rightarrow +\infty }} \frac{\mathbb{E}(\mathcal{D}(\mathbf{0},\mathcal{P}_n^d))}{n} .$$ 
Dhar ([Dha88]) obtained numerical upper bounds on $\mathbb{E}[\mathcal{D}(\mathbf{0},\mathcal{P}_n^d)]$ for small $n=1,2$ and valid for any $d$. This yields good upper bounds for $\mu_{\text{axis}}^d$.
Due to computer limitations, Dhar was not able to use his method in 1988 for $n>2$. The aim of this letter is to detail how to extend his computations, and how to compute efficiently new upper bounds for $\mu_{\text{axis}}^d$.

For all $n \in \mathbb{N}$ we note $\mathcal{J}_n^d$ the hyperplane $\mathcal{J}_n^d = \{ x_1 + x_2 + ... + x_d = \lfloor{} n\sqrt{d} \rfloor{} \}$. 
Observe that  $\mathcal{J}_n^d$ is chosen so that it is at the same Euclidean distance from $\mathbf{0}$ as $\mathcal{P}_n^d$. Therefore, if the Eden conjecture were true, one would have $\mathcal{D}(\mathbf{0},\mathcal{P}_n^d)= \mathcal{D}(\mathbf{0},J_n^d)+\mathrm{o}(n)$.
For the same reasons as before $\lim\limits_{n \rightarrow +\infty} \frac{\mathcal{D}(\mathbf{0},J_n^d)}{n}$ exists and thus we can define   $\mu_{diag}^d$ as $\mu_{diag}^d =\lim\limits_{n \rightarrow +\infty} \frac{\mathcal{D}(\mathbf{0},J_n^d)}{n}$. Couronné, Enriquez and Gerin (\citep{Gerin}) found an numerical lower bound on $\mu_{diag}^d$ (in fact, the same result also appears in a different form in \citep{dhar1986asymptotic}): 
$$\mu_{diag}^d \geq \frac{0.3313...}{\sqrt{d}}.$$

This means that a lower bound on the time of infection along the diagonal has been found. 
By combining his results with Dhar, \citep{Gerin} showed that  $\mu_{axis}^{35} < \mu_{diag}^{35}$, 
and thus that the limit shape of the infection $B^*$ is not an euclidean ball in dimension $35$.
In this letter we extend Dhar's method and use \citep{Gerin} lower bound to prove that $\mu_{axis}^{22} < \mu_{diag}^{22}$.

\begin{theorem} 
$\mu_{axis}^{22} < \mu_{diag}^{22}$. In particular the limiting shape $B^*$ is not an Euclidean ball, and the Eden conjecture is false in dimension 22.
\end{theorem}

\section{Dhar's strategy for $\mathcal{P}_1^d$}

\textbf{Idea} Since we will push Dhar's strategy a little further, we first detail the idea introduced in \citep{Dhar}. 
To compute an upper bound for $\mathbb{E}([\mathcal{D}(\mathbf{0},\mathcal{P}_n^d)])$, Dhar slightly modifies the model and considers a \emph{unidirectional} infection. This means that a site in $\mathcal{P}_i^d$ can only contaminate its neighbours in $\mathcal{P}_{i}^d$ and $\mathcal{P}_{i+1}^d$. We note $\tau_n^d$ the time of infection from $\mathbf{0}$ to the plan $\mathcal{P}_n^d$ by the \emph{unidirectional} infection. It is clear that this infection spreads more slowly than the original model, and therefore we obtain for every $n$
$$
\mu_{\text{axis}}^d \leq \frac{\mathbb{E}([\mathcal{D}(\mathbf{0},\mathcal{P}_n^d)])}{n} \leq  \frac{E[\tau_n^d]}{n}.
$$

From now on, we only consider the model of unidirectional infection in our computations.


\subsection{Notations}

\begin{figure}[h!]
\centering
\includegraphics[width = 0.4 \textwidth]{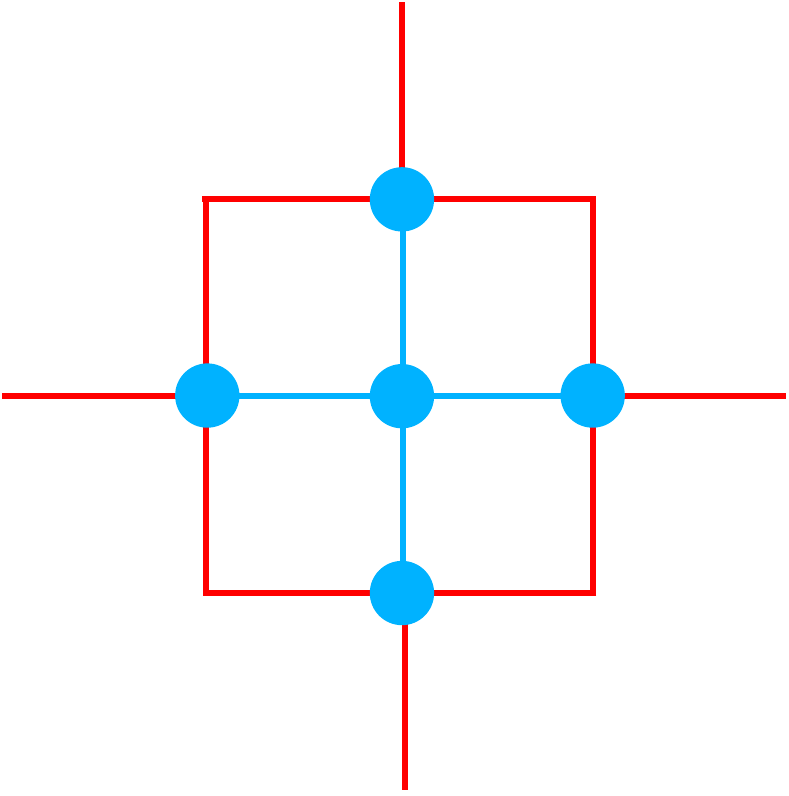}
\caption{The cluster $C$ (in blue) and the set of perimeter bounds $S$ (in red)}
\label{fig:schemaS}
\end{figure}

We consider a $d$-dimensional infection. Let $T(C)=\mathbb{E}(\tau_1^d | B_0 = C)$ be the expected waiting time before the infection reaches $\mathcal{P}_1^d$ starting from an infected cluster $C$ in $\mathcal{P}_0^d$ (and the other sites are healthy).
For a cluster $C \subset \mathcal{P}_0^d$ of $i$ sites, we define $S$ its set of perimeter edges (see Figure ~\ref{fig:schemaS}). As stated in \citep{Dhar} we have 

\begin{equation}
    \card{S} \geq \lceil 2(d-1)i^{\frac{d-2}{d-1}} \rceil := s_i.
\end{equation}


\begin{figure}[h!]
\centering
\includegraphics[width = 0.7 \textwidth]{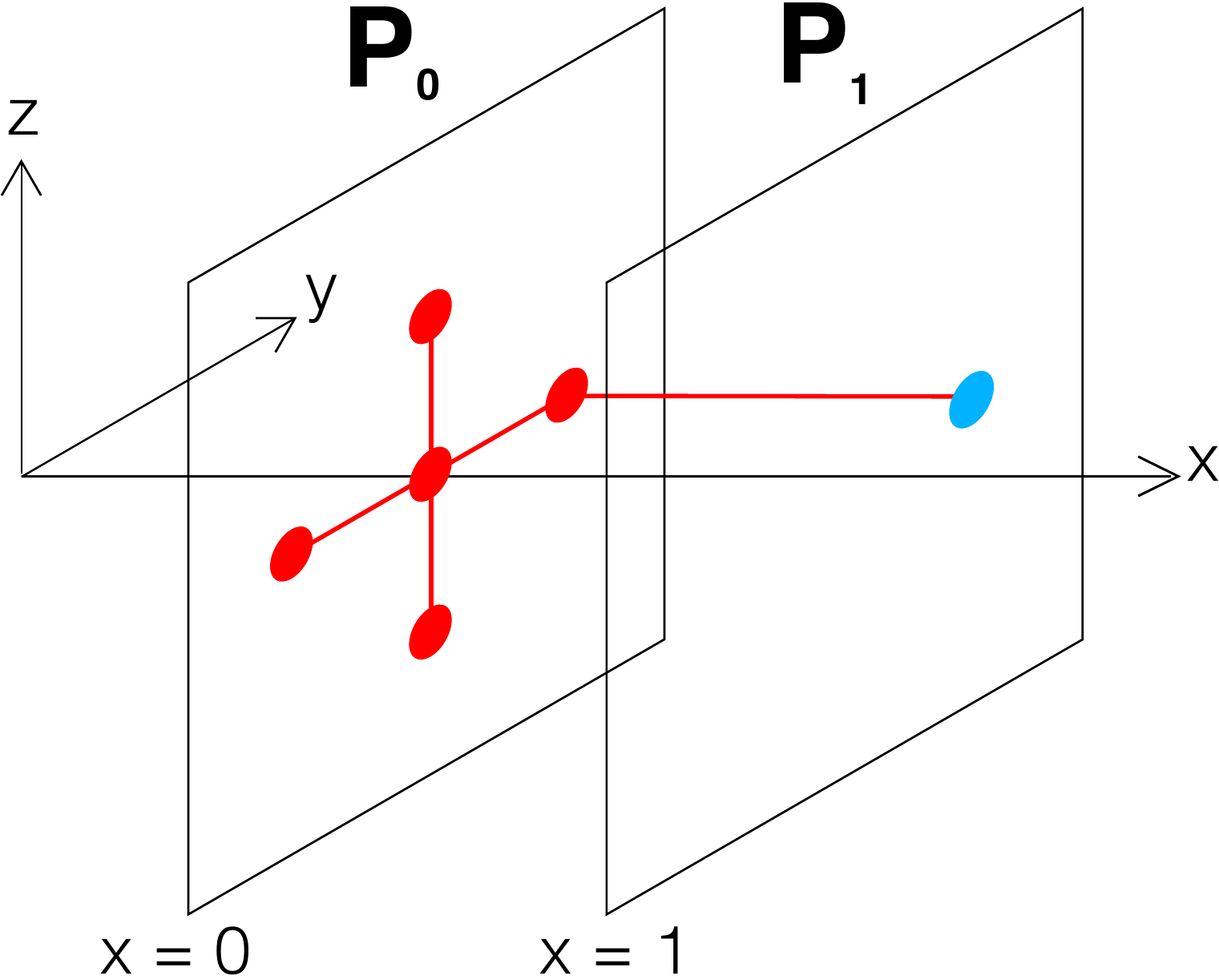}
\caption{Example of infection. In red the starting cluster $C$.}
\label{fig:schema2plan}
\end{figure}

For an edge $e=(x,y)$ such that $x\in C$ and $y\notin C$, we set $v^+(e)=y$ (\emph{i.e.} $v^+(e)$ is the endpoint of $e$ which is not in $C$).
We define $T_i = \displaystyle \max_{\substack{\card{C} = i}} T(C)$.




\subsection{Recursive inequality}


Let  $t_1$ be the time at which the first contamination occurs. At time $t_1$, a site in cluster $C$ contaminates either one of its neighbours in $\mathcal{P}_0^d$ or one of its neighbours in $\mathcal{P}_1^d$ (see figure~\ref{fig:schema2plan}). The total number of such "susceptible" sites is given by $|S|+i$, it follows that $t_1\stackrel{\text{def}}{=}\min \{z_1,\dots,z_{|S|+i}\}$ where $z_i$'s are i.i.d. passage times, \emph{i.e.} $t_1$ is distributed as an exponential r.v. with mean $1/(|S|+i)$. Using the same argument we have: 
\begin{equation}
    T_i \leq \dfrac{1}{i}.
    \label{eq:borneTi}
\end{equation}
At time $t_1^+$, the new infected site $x$ is uniformly distributed among the $|S|+i$ possibilities.
If $x\in \mathcal{P}_1^d$ then $\tau_1=t_1$. If $x\in \mathcal{P}_0^d$, then the infection goes on, starting from configuration $C'=C\cup\{x\}$. Because of the memoryless property of the exponential distribution, we have the Markov property
$$
T(C)-t_1\ |\ \{x \text{ is infected at time }t_1\}\ \stackrel{\text{def}}{=} T(C\cup\{x\}).
$$
Therefore we obtain

\begin{align*}
T(C) &= \underbrace{\mathbb{E}\left( t_1 \middle| B_0 = C \right)}_{\frac{1}{\card{S}+i} \text{ : min of $\card{S}+i$ exponential r.v}} +  \mathbb{E}\left(\tau_1^d - t_1  \middle| B_0 = C \right) & \\
&= \frac{1}{\card{S}+i} + \displaystyle \sum\limits_{\substack{e \text{ leaving } C}}{\mathbb{E}\left( \tau_1^d - t_1  \middle| B_0 = C, B_{t_1} = C' \right)\mathbb{P}\left(C' = C \cup v^+\{e\} \middle| B_0 = C \right)} \\
&= \frac{1}{\card{S}+i} + \displaystyle \sum\limits_{\substack{ e \text{ leaving } C \\ v^+\{e\} \in \mathcal{P}_{0}}}{\mathbb{E}\left( \tau_1^d - t_1  \middle| B_0 = C, B_{t_1} = C' \right)\mathbb{P}\left(C' = C \cup v^+\{e\} \middle| B_0 = C \right)}\\
&+ \displaystyle \sum\limits_{\substack{e \text{ edge leaving } C \\ v^+\{e\} \in \mathcal{P}_{1}}}{\underbrace{\mathbb{E}\left( \tau_1^d - t_1  \middle| B_0 = C, B_{t_1} = C' \right)}_{0 \text{ : } \tau_1^d = t_1  \text{ since we have reached } \mathcal{P}_1 }\mathbb{P}\left(C' = C \cup v^+\{e\} \middle| B_0 = C \right) }\\
&=  \frac{1}{\card{S}+i} + \displaystyle \sum\limits_{\substack{e \text{ edge leaving } C \\ v^+\{e\} \in \mathcal{P}_{0}}}{\underbrace{\mathbb{E}\left( \tau_1^d - t_1  \middle| B_0 = C, B_{t_1} = C' \right) }_{\mathbb{E}\left( \tau_1^d \middle| B_0 = C' \right)\text{ by Markov property}}\underbrace{\mathbb{P}\left(C' = C \cup v^+\{e\} \middle| B_0 = C \right)}_{\frac{1}{\card{S}+i} \text{ : choice of one edge among } \card{S}+i}} \\
&=  \frac{1}{\card{S}+i}\left( 1 + \displaystyle \sum_{e \text{ edge leaving } C}{T\left( C \cup v^+\{e\}\right)} \right).\\
\end{align*}

We have 

\begin{equation}
    T(C) \leq \dfrac{1+\card{S}T_{i+1}}{\card{S}+i}.
\end{equation}

The right-hand side is decreasing in $\card{S}$ and $|S| \leq s_i$ thus

\begin{equation}
    T_{i} \leq \dfrac{1+s_i T_{i+1}}{s_i+i}.
\end{equation}

which is inequality (8) in \citep{Dhar} (note a small misprint in Dhar's inequality (8)).
This recursive inequality and a rough bound on $T_n$ leads to a tight bound on $T_1$:
\begin{equation}
   \mathbb{E}(\tau_1^d) = T_1 \leq \frac{1+s_1 T_2}{1+s_1} \leq \frac{1+s_1 \frac{1+s_2 T_3}{2+s_2} }{1+s_1} \leq \frac{1+s_1 \frac{1+s_2 \frac{1+s_3 T_4}{3+s_3}}{2+s_2} }{1+s_1} \leq ...
\end{equation}









The numerical results for the bound on $\mathbb{E}(\tau_1^d)$ can be found in subsection \ref{results}, page \pageref{results}.

\section{An upper bound for $\mathbb{E}(\tau_2^d), \mathbb{E}(\tau_3^d), \mathbb{E}(\tau_4^d), \mathbb{E}(\tau_5^d) $}

\textbf{Idea} We consider the same unidirectional infection. We still use  $\mu_{axis}^d \leq \frac{\mathbb{E}(\mathbb{\mathcal{D}}(\mathbf{0},\mathcal{P}_{n}^d))}{n} \leq \frac{\tau_n^d}{n}$ but now we found bounds for $\tau_n^d, \ n \geq 1$ considering a starting cluster not just in $\mathcal{P}_0$ but in $\mathcal{P}_0^d, \dots, \mathcal{P}_{n-1}^d$.

\subsection{Upper bound on $\mathbb{E}(\tau_2^d)$}

We consider an modified infection starting with starting clusters $C_0 \subset \mathcal{P}_0^d$ with $i$ sites and $C_1 \subset \mathcal{P}_1$ with $j$ sites. For $k \in \{0, 1\}$ we note $S_k$ the set of perimeter edges connecting the cluster $C_k$ to the rest of the plan $\mathcal{P}_k^d$ see Figure \ref{fig:schemaS}.
Similarly to the previous section, we define $T(C_0,C_1)=\mathbb{E}(\tau_2^d | B_0 = C_0 \cup C_1 )$ and $$T_{i,j} = \displaystyle \max_{\substack{\card{C_0} = i \\ \card{C_1} = j }} T(C_0,C_1)$$

\begin{figure}[h!]
\centering
\includegraphics[width = 0.5 \textwidth]{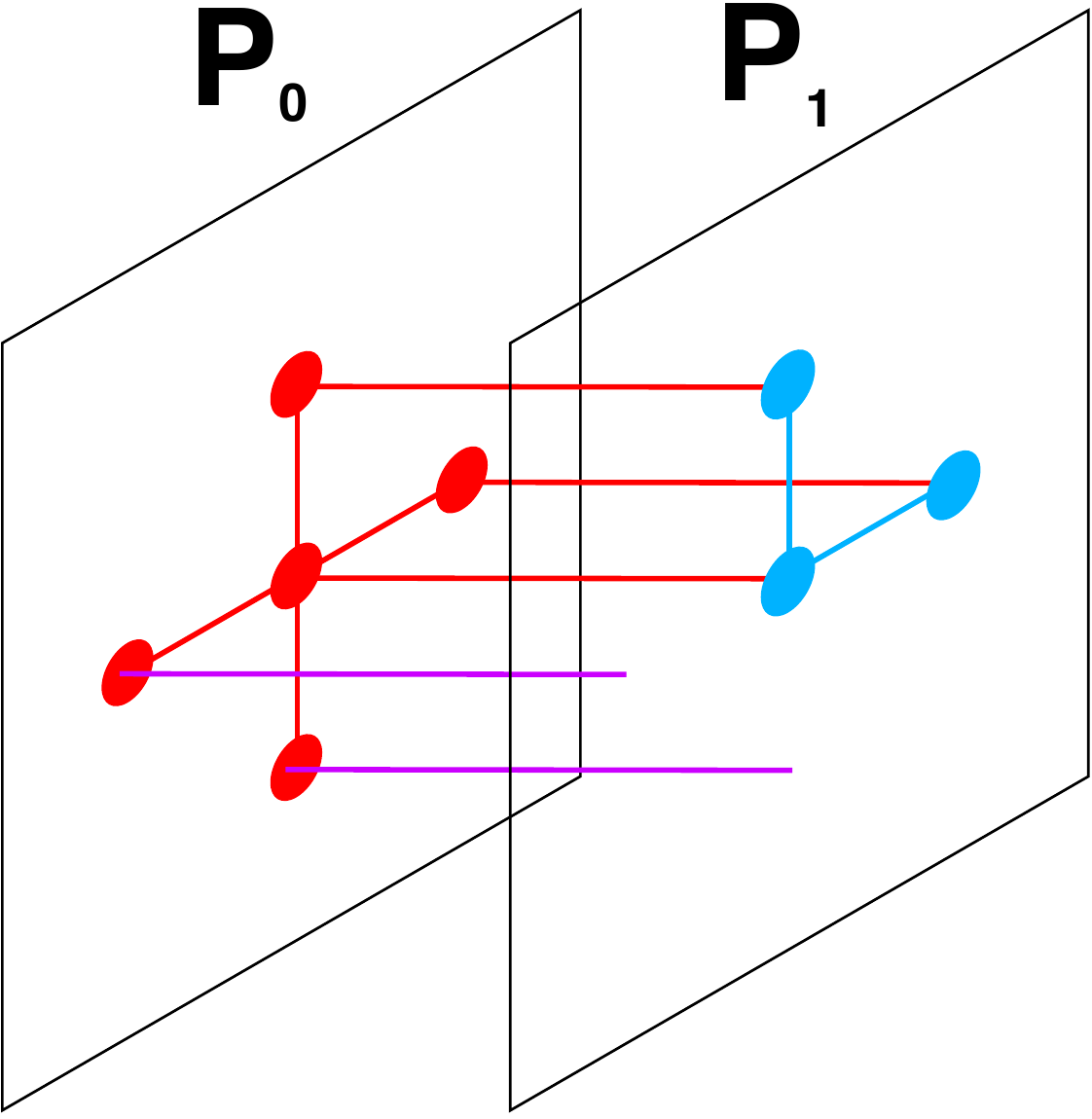}
\caption{The starting clusters $C_0$ (in red) and $C_1$ (in blue) and the edges of B (in purple).}
\label{fig:schemaB}
\end{figure}

Le $B$ be the set of edges from $C_0$ to healthy sites in $\mathcal{P}_1^d$ (see figure~\ref{fig:schemaB}). Clearly, $|B|\leq i$. On the other hand, $B$ is minimal when $C_0$ is exactly in front of $C_1$: therefore $|B| \geq |i - j|_+$.


Using the same method that in the previous section we have 

\begin{equation}
    T(C_0,C_1) = \frac{\left( 1 
    + \displaystyle \sum_{e \text{ edge leaving } C_0}{T\left(C_0' = C_0 \cup v^+\{e\}, C_1\right)} + \displaystyle \sum_{e \text{ edge leaving } C_1}{T\left(C_0, C_1' = C_1 \cup v^+\{e\}\right)}\right)}{\card{S_0}+\card{S_1}+\card{B}+j}.
\end{equation}

This leads to 
\begin{equation}
    T(C_0, C_1) \leq \displaystyle  \frac{ 1+ \card{S_0} T_{i+1,j} + \left(\card{S_1}+\card{B}\right)T_{i,j+1}} {\card{S_0}+\card{S_1}+\card{B}+j}.
    \label{eq:recuTau2_0}
\end{equation}

Like in the previous section the right hand-side of the inequality is decreasing in the variables $\card{S_0}$ and $\card{S_1}+\card{B}$ (the proof is in annex). Morover $\card{S_0} \geq s_i$ and $\card{S_1}+\card{B} \geq s_j + |i-j|_+ $. This leads to 

\begin{equation}
    T_(C_0, C_1) \leq \displaystyle \frac{ 1+ s_i T_{i+1,j} + \left(s_j + |i-j|_+ \right)T_{i,j+1}} {s_i+s_j+|i-j|_+ +j}
    \label{eq:recuTau2_1}.
\end{equation}

And by taking the maximum of $T_(C_0, C_1)$ in the set $\{ C_0, C_1 \text{ such that } |C_0|=i, |C_1|=j \}$ it leads to

\begin{equation}
    T_{i,j} \leq \displaystyle \frac{ 1+ s_i T_{i+1,j} + \left(s_j + |i-j|_+ \right)T_{i,j+1}} {s_i+s_j+|i-j|_+ +j}
    \label{eq:recuTau2}.
\end{equation}
 


With the recursion inequality \ref{eq:recuTau2}. It is now possible to recursively compute an upper bound on $T_{1,0} = \mathbb{E}(\tau_2)$.

\subsection{Boundary conditions}

The computation is also done backwards. For initialization, we use the upper bounds for $T_i$ that we obtained in the previous section.

\begin{figure}[h!]
\centering
\includegraphics[width = 0.5 \textwidth]{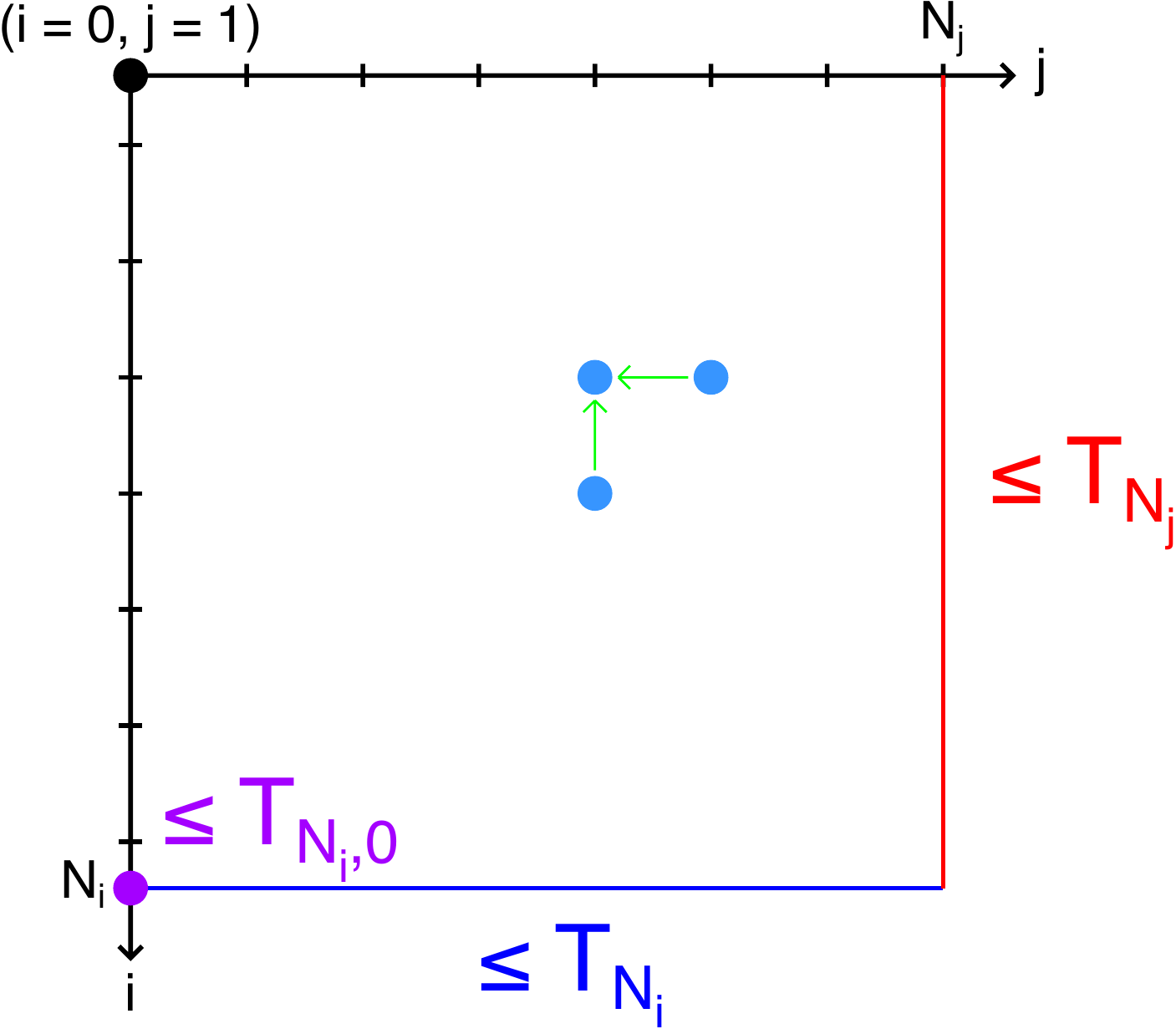}
\caption{Filling process of the array $T_{i,j}$ while computing the bound on $\mathbb{E}(\tau_2)$.}
\label{fig:remplissage2D}
\end{figure}

As shown in Fig. \ref{fig:remplissage2D} we roughly bound borders and then appply the backward inequation \ref{eq:recuTau2}:

\[
\left\{
\begin{array}{lclc}
     T_{i,jMax} & \leq & T_{jMax} & \forall i \in \llbracket 1, iMax \rrbracket . \\
     T_{iMax, j} & \leq & T_{j} & \forall j \in \llbracket 1, jMax \rrbracket.   \\
     T_{iMax, 0} & \leq & \min \left\{ \gamma^1, \dots, \gamma^{iMax} \right\} & \gamma^k \overset{iid}{\sim} \text{gamma(2,1).}
\end{array}
\right.
\]

The last inequality is obtained by considering the $iMax$ disjoint paths $x \rightarrow x+(1,0,0) \rightarrow x+(2,0,0)$ for $x$ in $C_0$. The length of each such path is the sum of two independent exponential r.v., \emph{i.e.} a Gamma(2,1). The numerical results can be found in subsection \ref{results}, page \pageref{results}.

\subsection{Expansion to $\mathbb{E}(\tau_n^d)$}

Using the same technique, we can compute the upper bound $\mathbb{E}(\tau_n^d)$ for all $n$.We have to dynamically fill a n-dimensionnal array which remains reasonnably time consuming for $n\leq 5$. We still get initialising bounds with the previous calculation and the straight forward paths to $\mathcal{P}_n$. 

For $\mathbb{E} (\tau_3)$ we obtain a dynamic equation of level 3 with the following bounding conditions:


\[\left\{
\begin{array}{lclc}

    T_{iMax,j, k} & \leq & T_{j, k} & \forall (j, k) \in \llbracket 1, jMax \rrbracket \times \llbracket 0, kMax \rrbracket\\
    T_{iMax,0, k} & = & 0 & \forall k \in  \llbracket 1, kMax \rrbracket\\
    T_{i,0, kMax} & = & 0 & \forall i \in  \llbracket 1, iMax \rrbracket\\
    T_{iMax,0, k} & = & 0 & \forall k \in  \llbracket 1, kMax \rrbracket\\
     T_{i,jMax, k} & \leq & T_{jMax, k} & \forall (i, k) \in \llbracket 1, iMax \rrbracket \times \llbracket 0, kMax \rrbracket\\
     T_{i,j, kMax} & \leq & T_{j, kMax} & \forall (i,j) \in \llbracket 1, iMax \rrbracket \llbracket 1, jMax \rrbracket   \\

     T_{iMax, 0, 0} & \leq & \min \left\{ \gamma^1, \dots, \gamma^{iMax} \right\} & \gamma^k \overset{iid}{\sim} \text{gamma(3,1)}.
\end{array}
\right.\]

For the numerical computations of $\mathbb{E}(\tau_4^d)$ and $\mathbb{E}(\tau_5^d)$, the method is the same.
We did not go any further in the calculus because filling the hypercube began to be very time consuming.  The numerical results can be found in subsection \ref{results}, page \pageref{results}.

\section{Implementation and numerical results}

\subsection{Numerical Tricks}

\textbf{Sensitivity to boundary data.} To get the upper bound of $\mathbb{E}(\tau_n)$, one can see from the dynamic equations that we use the $T_{i_1, \ldots, i_{n-1}}$ previously computed for $\mathbb{E}(\tau_{n-1})$ to initialize the boundary of the hypercube. We saw numerically that the resulting upper bound is very sensitive to this initialisation. With a better precision of this boundary data, we get far better upper bounds. 


\textbf{Choice of parameters $\mathrm{iMax},\mathrm{jMax}$, $\mathrm{kMax}$.} We used this property to efficiently chose the size of the box. It's more efficient to put the computational effort on boundary data than on the recursion. For instance, imagine we want to compute the upper bound on $T_{1,0,0}$. We have recursion with order 3, and limit conditions involving the calculus of $T_{i,j}$. To compute the bound on $T_{1,0,0}$ we may want to build a $1000 \times 1000 \times 1000$ box to fill it.
But we can obtain faster results by filling a $1000 \times 100 \times 100$ box, and using a $10000 \times 1000$ box to compute the limits-condition $T_{i,j}$.


\subsection{Numerical results}

\label{results} 

We did the implementation for the computation of the bounds on $\mathbb{E}(\tau_1^d), ..., \mathbb{E}(\tau_5^d)$ in \texttt{C++} \footnotemark{} \footnotetext{\url{https://github.com/QB3/Dimension-improvement-in-Dhar-s-refutation-of-the-Eden-conjecture}}. Using this code we obtained the following results:


\begin{tabular}{|c|c|c|c|c|c|c|}
\hline
Dimension d &  $ \mu_{diag}^d \geq$  & $\mathbb{E}(\tau_1^d) \leq $  & $\frac{\mathbb{E}(\tau_2^d)}{2} \leq $ & $\frac{\mathbb{E}(\tau_3^d)}{3} \leq $ & $\frac{\mathbb{E}(\tau_4^d)}{4} \leq $  & $\frac{\mathbb{E}(\tau_5^d)}{5} \leq$\\
\hline
2 & 0.2343 &  0.5973 & 0.5560 &0.5341  & 0.5211 &0.5126 \\
3 & 0.1913 & 0.4400 & 0.4011 & 0.3813 & 0.3697 & 0.3622 \\
4 & 0.1657 & 0.3527 & 0.3177 & 0.3004 & 0.2903 &0.2839  \\
5 & 0.1482 & 0.2976 & 0.2662 & 0.2507  & 0.2419 &0.2362  \\
22 & \textcolor{red}{0.0706} & 0.0933 & 0.0812 & 0.0758 & 0.0730 & \textcolor{red}{0.0699}  \\
25 & \textcolor{red}{0.0663} & 0.0842 & 0.0732 & 0.0684 & \textcolor{red}{0.0660} &0.0644  \\
30 & \textcolor{red}{0.0605} & 0.0727 & 0.0631 & \textcolor{red}{0.0592} & 0.0572 &0.0560  \\
35 & \textcolor{red}{0.0560} & 0.0642 & \textcolor{red}{0.0556} & 0.0524  & 0.0507 &0.0497  \\
\hline
\end{tabular}

\subsection{Future work}

We focused on $\mu^d_{axis}$ and did not work at all on $\mu^d_{diag}$. But in order to prove the conjecture wrong in dimension 2 bounds on $\mu^2_{axis}$ and $\mu^2_{diag}$ need to be improved. Indeed Monte-Carlo simulations (\citep{alm2015first}) give $\mu^2_{axis}  \simeq 0.404$ and $\mu^2_{diag} \simeq 0.409$. To improve $\mu^d_{axis}$ Dhar suggested us a trick: if $\mathbb{E}(T_{n+1}-T_n) \leq \mathbb{E}(T_{n}-T_{n-1})$ then it can easily be shown that $\mu^2_{axis}  \leq \mathbb{E}(T_{n}-T_{n-1})$. This lead to a slightly different recursion inequality on a bound for $\mathbb{E}(T_{n}-T_{n-1})$ which allowed us to catch the dimension $d=16$. Unfortunately neither Dhar or us finaly managed to find a correct proof of $\mathbb{E}(T_{n+1}-T_n) \leq \mathbb{E}(T_{n}-T_{n-1})$.

\section*{Acknowledgements}

We would like to thank Lucas Gerin for introducing us to the subject. He kindly guided us through this project, and took time to answer all our questions. We warmly thank Deepak Dhar for his precious help and the huge amount of time he spent for us.

\appendix
\section{Annex}

\textbf{Goal: } prove that \begin{equation}
    \displaystyle  \frac{ 1+ \card{S_0} T_{i+1,j} + \left(\card{S_1}+\card{B}\right)T_{i,j+1}} {\card{S_0}+\card{S_1}+\card{B}+j}.
    \label{eq:recuTau2_0} \end{equation} is decreasing in $\card{S_0}$ and $\card{S_1}+\card{B}$ 

We proved that: $T(C_0, C_1) \leq \displaystyle  \frac{ 1+ \card{S_0} T_{i+1,j} + \left(\card{S_1}+\card{B}\right)T_{i,j+1}} {\card{S_0}+\card{S_1}+\card{B}+j}. $ Thus

$ \underset{C_0, C_1\text{ such that }|C_0|=i\text{ and }C_1=j}{\max} T(C_0, C_1) \leq \displaystyle  \underset{C_0, C_1\text{ such that }|C_0|=i\text{ and }C_1=j}{\max} \frac{ 1+ \card{S_0} T_{i+1,j} + \left(\card{S_1}+\card{B}\right)T_{i,j+1}} {\card{S_0}+\card{S_1}+\card{B}+j}. $ 

The maximization is over a finite set, this means that there exists $S_{0}^*$, $S_1^*$ and $B^*$ such that 
$ \frac{ 1+ \card{S_0^*} T_{i+1,j} + \left(\card{S_1^*}+\card{B^*}\right)T_{i,j+1}} {\card{S_0^*}+\card{S_1^*}+\card{B^*}+j} = \underset{C_0, C_1\text{ such that }|C_0|=i\text{ and }C_1=j}{\max} \frac{ 1+ \card{S_0} T_{i+1,j} + \left(\card{S_1}+\card{B}\right)T_{i,j+1}} {\card{S_0}+\card{S_1}+\card{B}+j}.$

We can rewrite the inequality as: $T_{i,j} \leq \frac{ 1+ \card{S_0^*} T_{i+1,j} + \left(\card{S_1^*}+\card{B^*}\right)T_{i,j+1}} {\card{S_0^*}+\card{S_1^*}+\card{B^*}+j}.$ 

It follows $T_{i,j} \leq \frac{ 1+ (\card{S_0^* + \card{S_1^*}+\card{B^*}+j - \card{S_1^*}- \card{B^*} - j})  T_{i+1,j} + \left(\card{S_1^*}+\card{B^*}\right)T_{i,j+1}} {\card{S_0^*}+\card{S_1^*}+\card{B^*}+j}.$ And then 

$T_{i,j} \leq T_{i+1,j} + \frac{ 1- ( \card{S_1^*} + \card{B^*} +j)  T_{i+1,j} + \left(\card{S_1^*}+\card{B^*}\right)T_{i,j+1}} {\card{S_0^*}+\card{S_1^*}+\card{B^*}+j}.$

However we know that $T_{i+1,j} \leq T_{i,j}$ this means that $\frac{ 1- ( \card{S_1^*} + \card{B^*} +j)  T_{i+1,j} + \left(\card{S_1^*}+\card{B^*}\right)T_{i,j+1}} {\card{S_0^*}+\card{S_1^*}+\card{B^*}+j} \geq 0$

Let now study the function $f :  x \mapsto \frac{ 1- (\card{S_1^*} + \card{B^*} +j)  T_{i+1,j} + (\card{S_1^*} + \card{B^*} )T_{i,j+1}} {x +\card{S_1^*} + \card{B^*}+j}.$

$f$ is monotone, tends to $0$ when $x$ goes to infinity moreover we proved that there exists $S_0^*$ such that $f(S_0^*) \geq 0$. This proves that f is decreasing. We can replace $S_0^*$ by a lower bound on it: $s_i$. This leads to 

$T_{i,j} \leq T_{i+1,j} + \frac{ 1- ( \card{S_1^*} + \card{B^*} +j)  T_{i+1,j} + \left(\card{S_1^*}+\card{B^*}\right)T_{i,j+1}} {s_i+\card{S_1^*}+\card{B^*}+j}.$

We can use exactly the same trick to prove that $\card{S_1^*} + \card{B^*}$ can be replaced by a lower bound: $s_j + |i-j|^+$

\bibliographystyle{alpha}
\bibliography{Eden}

\begin{thebibliography}{ADH15}

\bibitem[AD15]{alm2015first}
Sven~Erick Alm and Maria Deijfen.
\newblock First passage percolation on {$\mathbb{Z}^{2}$} : A simulation study.
\newblock {\em Journal of statistical physics}, 161(3):657--678, 2015.

\bibitem[ADH15]{refPerco}
Antonio Auffinger, Michael Damron, and Jack Hanson.
\newblock 50 years of first passage percolation.
\newblock {\em arXiv preprint arXiv:1511.03262}, 2015.

\bibitem[CEG11]{Gerin}
Olivier Couronn{\'e}, Nathana{\"e}l Enriquez, and Lucas Gerin.
\newblock Construction of a short path in high-dimensional first passage
  percolation.
\newblock {\em Electronic Communications in Probability}, 16:22--28, 2011.

\bibitem[Dha86]{dhar1986asymptotic}
Deepak Dhar.
\newblock Asymptotic shape of eden clusters.
\newblock In {\em On growth and form}, pages 288--292. Springer, 1986.

\bibitem[Dha88]{Dhar}
Deepak Dhar.
\newblock First passage percolation in many dimensions.
\newblock {\em Physics Letters A}, 130(4-5):308--310, 1988.

\bibitem[HW65]{hammersley1965first}
John~M Hammersley and DJA Welsh.
\newblock First-passage percolation, subadditive processes, stochastic
  networks, and generalized renewal theory.
\newblock In {\em Proc. Internat. Res. Semin., Statist. Lab., Univ. California,
  Berkeley}, pages 61--110. Springer, 1965.

\end{thebibliography}







\nocite{*}

\end{document}